
\documentstyle{amsppt}
\magnification \magstep1
\parskip 11pt
\parindent .3in
\pagewidth{5.3in}
\pageheight{7.2 in}

%
%
\def \bl{\vskip 11pt}

\def \ni{\noindent}
\def \P{\bold{P}}

\def \Z{\bold{Z}}

\def \C{\bold{C}}
\def \O{\Cal{O}}
\def \U{\Cal{U}}
\def \T{\Bbb T}
\def \W{\Cal W}
\def \H{\Cal H}

\def \M{\Cal {M}}
\def \B{\Cal B}
\def \E{\Cal{E}}
\def \N{\Bbb N}
\def \L{\Cal L}
\def \K{\Bbb K}
\def \nilp{\text{Nilp}}
\def \overbar{\overline}
\def \lra{\longrightarrow}
%

\centerline{\bf A NON-LINEAR DEFORMATION OF THE HITCHIN DYNAMICAL SYSTEM}
\bl
\bl
\centerline{\smc Ron Donagi \footnote{Partially supported by  NSA Grant
MDA904-92-H3047 and NSF Grant DMS 95-03249 }}
\centerline{\smc Lawrence Ein \footnote{Partially supported by NSF Grant DMS
93-02512}}
\centerline{\smc Robert Lazarsfeld \footnote{Partially supported by NSF Grant
DMS 94-00815}}
\bl
\bl
\ni{\bf Introduction.}  The purpose of this note is to relate two
constructions that have attracted considerable attention among algebraic
geometers. On the one hand, Hitchin \cite{H} has introduced a beautiful
completely integrable dynamical system living on the cotangent bundle
$\T^*\U_C$ of the moduli space $\U_C$ of stable vector
bundles of given rank and degree on an algebraic curve $C$.  On the other
hand, Mukai \cite{M} has studied the space $\M$ parametrizing stable sheaves
of given numerical type on a $K3$ surface $S$, proving in particular that
$\M$ is a symplectic variety. We will argue here that it is profitable to view
(certain cases of) the Mukai contsruction as  a deformation of the Hitchin
system. The idea is that if $C$ is a hyperplane section of a $K3$-surface
$S$, then $S$ can be deformed to the cone $K$ over the canonical embedding of
$C$. Our deformation of $\M$ to the Hitchin system (given in \S1) is based on
exactly this deformation of $C \subset S$ to its normal cone.

 We apply this deformation in \S 2 to study the ``nilpotent cone" Nilp$(\M)$
in the Mukai system, which parametrizes certain sheaves on a  non-reduced
curve on  $S$. Drawing on Laumon's analysis \cite{L} of the corresponding
locus in $\T^*\U_C$, we show that Nilp$(\M)$ is a Lagrangian subvariety of
$\M$ provided at least that certain numerical restrictions are satisfied.
Finally, in \S3, we sketch an elementary result describing in detail the
nilpotent cones of the two systems  in the rank two case.  In  Hitchin's
setting, this consists of a union of vector bundles over various symmetric
products of the base curve
$C$. In the Mukai space, the corresponding locus consists of {\it affine}
bundles having the same underlying linear structure as in the Hitchin case.

Consideration of the Mukai and Hitchin systems as moduli spaces of sheaves
on a fixed surface suggests some general questions as to when such a
moduli space has a symplectic, Poisson, or completely integrable structure.
One beautiful answer has recently been discovered by Markman \cite{Mkm}:  the
moduli space of ``Lagrangian" sheaves on a symplectic variety is a
completely integrable system, fibered by the support map. We discuss this and
some open questions in the final section.

We work throughout over the complex numbers $\C$.  If $E$ is a vector bundle
on a variety $X$, viewed say as a locally free sheaf, we denote by ${\Bbb E}
\lra X$ the total space of $E$, considered as an algebraic variety mapping to
$X$.

\bl
\ni {\bf \S 1. The systems of Hitchin and Mukai, and the Deformation.}  We
start by briefly recalling Hitchin's construction. Let $C$ be a fixed smooth
projective curve of genus $g \ge 2$, let $K_C$ be the canonical bundle of
$C$, and denote by
$\U = \U_C(n,e)$  the moduli space parametrizing stable rank $n$ vector
bundles on $C$, of degree $e$. Then $\U$ is a smooth variety, of dimension
$\tilde g =_{\text{def}}n^2(g-1) + 1$. In a well-known manner, one can identify
the
cotangent bundle $\T^* \U$ with the set of all pairs consisting of a stable
rank $n$ vector bundle (up to isomorphism) plus a $K_C$-valued endomorphism
of that bundle:
$$
\T^* \U = \left \{ (E, \phi) \mid \phi : E \lra E \otimes K_C \right \}. $$
Let $B = \Gamma(C, K_C) \oplus \Gamma(C, K_C^2) \oplus \dots \oplus \Gamma(C,
K_C^n)$. Then one has a natural morphism of varieties:
$$ H : \  \T^* \U  \lra     B = \C^{\tilde g}, \tag 1.1
$$ given by $\ h(E, \phi ) = P_\phi$, where $$P_\phi = (-\text{tr } \phi,
\text{tr } \Lambda^2 \phi, \dots, (-1)^n\text{ det }\phi) \in B$$ is the
``characteristic polynomial" of
$\phi$. Hitchin proves that $H$ is an algebraically completely integrable
Hamiltonian system with respect to the natural symplectic structure on
$\T^* \U$. It follows that if $P \in B$ is a general point, then $H^{-1}(P)$
can be realized as a Zariski-open subset of an abelian variety $J_P$. These
abelian varieties can also be explained via the spectral construction.
In fact, an element $P \in B$, viewed as a characteristic polynomial,
determines a {\it spectral curve}
$$ D_P \subset {\Bbb K}_C $$  mapping $n$-to-one onto
$C$. Fibre-wise over $x \in C$,  $D_P$ parametrizes  the  eigenvalues of
$\phi(x)$. Then for general $P \in W$, $J_P = Jac(D_P)$. We refer to
\cite{BNR} for a discussion of spectral curves, and for an algebro-geometric
approach to Hitchin's results.

Turning to Mukai's work, let $S$ be a smooth projective polarized $K3$
surface. Mukai's beautiful observation \cite{M} is that the moduli spaces
parametrizing simple sheaves on $S$ are smooth, and carry natural symplectic
structures. We will be concerned with the following special case. Fix a
smooth curve $C \subset S$ of genus $g \ge 2$, and suppose for simplicity
that $C$ is a very ample divisor on $S$. Given integers $k , n \in \Z$,
with $n \ge 1$,  denote by  $$\overline \M = \M_{|nC|}^k(S)$$  the moduli
space parametrizing pairs consisting of a divisor $$D  \in \overline B
=_{\text{def}} |nC|,$$ plus a semi-stable coherent sheaf $\E$ supported on
$D$ having the same Hilbert polynomial as a line bundle of degree $k$ on a
smooth member
$D^\prime\in |nC|$. We refer to \cite{S}, \S1, for the precise
definition of semi-stability, and for the proof that $\overline \M$ exists.
\footnote "*"{Note  that one the requirements for semi-stability is that
$\E$ have ``pure dimension one", i.e. that it be free of embedded
components.} We will say that such a sheaf $\E$ has {\it numerical rank one
and degree} $k$. For example, if $E$ is a semi-stable vector bundle of rank
$n$ and degree $k + (n^2 -n)(1-g)$ on $C$, then we may view $E$ as a
coherent sheaf on  the non-reduced curve $R = nC \subset \overline S$, and
the pair $(E, R)$ determines a element of $\overline \M$. Since a
semi-stable sheaf is simple, it follows from  \cite{M} that  the open subset
$\M \subset \overline \M$ consisting of stable sheaves is a symplectic
manifold. As in \cite{LeP}, \S 2.3, taking schematic supports via Fitting
ideals defines a morphism
$$
\overline M : \overline \M \lra \overline B = \P^{\tilde g}, \tag 1.2 $$
where as above $\tilde g =  \ \text{dim} \ |nC| =n^2(g-1) + 1$.

Just as in the Hitchin system, one has
\proclaim{\bf Lemma 1.3} For a smooth curve $D \in |nC|$, the
fibre $\overline M^{-1}([D]) = Pic^k(D)$ is a Lagrangian submanifold
of $\M$.  \endproclaim
\demo{Proof} One computes that $D$ has genus $\tilde g$, and hence  dim $\M
= 2 \tilde g = 2 \text{ dim }\overline M^{-1}(D)$ along $\overline
M^{-1}(D)$. So it is enough to show that the symplectic form vanishes on the
fibres.   Let $L$ be a line bundle of degree
$k$ on $D$. Then the corresponding point $[L] \in \overline \M$ lies in the
stable locus
$\M \subset \overline \M$. The tangent space $T_{[L]}\M$ to $\M$ at $[L]$ is
identified with Ext$^1_S(L,L)$, whereas the tangent space $T_{[L]}Pic(D)$ to
the fibre is given by the subspace  Ext$^1_D(L, L) = \text{ Ext}^1_D(\O_D,
\O_D)$. It follows from Mukai's construction \cite{M} that the restriction
of the symplectic form to $T_{[L]}Pic(D)$ factors through the
cup-product map
$$\text{ Ext}^1_D(L, L) \times \text{ Ext}^1_D(L, L) \lra \text{ Ext}^2_D(L,
L) .$$ But $\text{ Ext}^2_D(L, L) = 0$, and the lemma follows. \qed \enddemo

\demo {\bf Remark 1.4} In both the Hitchin and Mukai systems, the Jacobians
$J(D)$ of smooth curves $D \in |nC|$ are realized as Lagrangian submanifolds
of a symplectic manifold. It is shown in \cite{DM1} and \cite{DM2}  that such a
situation is
determined by a family of cubic tensors $c_D \in Sym^3 H^0(D, K_D)^*$. In
each setting, the cubic at $D$ is given by the extension class
$$e_D \in \text{ Ext}^1_D(K_D, T_D) = H^0(D, K_D^3)^*$$ of the normal bundle
sequence
$$
0 \lra T_D \lra T_F|D \lra N_{D/F} \lra 0,$$
where $F = \overline \K_C$ or $F = S$ as the case may be.

We now wish to show that one can degenerate the Mukai system (1.2) to the
Hitchin system. Note to begin with that Hitchin's system admits a
natural compactification, as follows. Let
$\overline \K_C$ denote the one-point compactification of $\K_C$, i.e. the
projective completion $\P(K_C \oplus \O_C)$ with the section at infinity
blown down to a point. (So if $C$ is non-hyperelliptic, $\overline \K_C$ is
the cone over the canonical embedding of $C$.) View $C \subset \K_C
\subset \overline \K_C$ as embedded by the zero-section. Then the spectral
curves  $D_P \subset \overline \K_C$ lie in the linear series $\overline B
=_{\text{def}}| nC | = \P^{\tilde g}$, and one may identify the affine space
$B$ as the subset of $\overline B $ formed by those curves not
passing through the vertex of $\overline \K_C$. Given an integer $k$, we
consider pairs consisting of a curve $D \in |nC|$ plus a semi-stable
coherent sheaf $\E$ supported on $D$ having the same Hilbert polynomial as a
line bundle of degree $k$ on a smooth member $D^\prime \in |nC|$. We will say
as above that such a sheaf $\E$ has  numerical rank one and degree  $k$.
According to Simpson \cite{S}, there is a projective moduli space
$$\overline \H = \overline \H ^k_{|nC|}(\overline \K_C)$$
parametrizing isomorphism classes of all such, where naturally we view
$\overline \K_C$ as polarized by the ample divisor $C$.  Taking supports
defines as before  a morphism
$$\overline H: \overline \H \lra \overline B. \tag 1.5$$
Via the spectral construction $\T^* \U_C$ is realized as an open subset of
$\overline H ^{-1}(B)$, with $\overline H$ restricting to the Hitchin map.
We denote by $\H \subset \overline \H$ the open subset parametrizing stable
sheaves.

Now consider  a $K3$ surface $S \subset \P^g$ containing a
smooth curve $C \subset S$ of genus $g$ as a hyperplane section. The point is
simply to exploit the elementary fact that one can degenerate $S$ to
the cone $\overline \K_C$ over $C$. Specifically, let $X_0 \subset \P^{g+1}$
be the cone over $S$, and denote by $X \lra X_0$  the blowing-up of $X_0$
along $C$. The pencil of hyperplanes in $\P^{g+1}$ passing through $C$ gives
rise to a mapping
$$f : X \lra \P^1 .$$  There is a distinguished point $O \in \P^1$ such that
$$f^{-1}(O) = \overline \K_C;$$  for all other points $O \ne t \in \P^1$,
$f^{-1}(t) \cong S$. The map $f$ is thus a deformation of $C \subset S$
to the cone over $C$. It was used for example in \cite{P}, and is
closely related to the well-known deformation (cf. \cite{F}) of $C
\subset S$ to its normal cone. (The latter is essentially the
blow-up of $X$ at the vertex of $\overline \K_C \subset X$.)  Note that
$C$ embeds naturally as an ample divisor in each of the fibres $X_t$ of $f$,
and that $X$ carries a polarization whose restriction to each fibre is
$\O_{X_t}(C)$.

Now recall that Simpson \cite{S} constructs moduli spaces of sheaves in a
relative setting. Let
$$\overline \W = \overline \W^k_{|nC|} \lra \P^1$$
be the moduli space, projective over $\P^1$, parametrizing semistable
sheaves of numerical rank one and degree $k$ contained in the fibres of
$f$.  Equivalently, $\overline \W$ may be described as the moduli space
parametrizing semi-stable sheaves of the given numerical type on $X$ which
are supported in a fibre of $f$. Denote by
$$\overline \B \lra  \P^1$$ the $\P^{\tilde g}$-bundle whose fibre over $t
\in \P^1$ is the linear series $| \O_{X_t}(nC)|$. Then Le Potier's
construction in \cite{LeP} globalizes to define a support map
$$
\overline W : \overline \W \lra \overline \B \tag 1.6$$
of schemes over $\P^1$. Given $t \in \P^1$, let $$\overline \W_t \lra
\overline \B_t$$ denote the the fibres of (1.6) over $t$. Then for $O \ne t
\in \P^1$  one has $\overline \W_t \cong \overline \M$, whereas $\overline
\W_O \cong \overline \H$. Thus (1.6) defines the required degeneration
of the Mukai to the Hitchin system.

\bl
\ni {\bf \S 2. The Nilpotent Cone in Mukai Space.} In the ``classical"
Hitchin system $H : \T^* \U_C \lra B$, the {\it nilpotent cone} consists of
all pairs $(E, \phi)$ such that the endomorphism $\phi : E \lra E \otimes
K_C$ is nilpotent. A basic theorem of Laumon \cite{L} states that this cone
is a Lagrangian subvariety of $\T^* \U_C$. In this section, we use the
deformation (1.6) to deduce the corresponding statement for the Mukai
system under certain numerical restrictions.

The first point is to define the nilpotent cones in the present setting. To
this end, note that if $\phi : E \lra E \otimes K_C$  is a nilpotent
endomorphism of a semi-stable bundle $E$ on $C$, then the corresponding
spectral curve is the non-reduced scheme $R = nC \subset \K_C$, i.e.  $R$ is
the
$(n-1)^{\text{st}}$ infinitesimal neighborhood of the zero section $C
\subset \K_C$. Therefore we define the nilpotent
variety $\text{Nilp}(\overline \H)  \subset \overline \H = \overline \H
^k_{|nC|}(\overline \K_C)$ to be the scheme
$$\text{Nilp}(\overline \H) = \overline H^{-1}([R])$$
parametrizing  all semi-stable sheaves of numerical rank one and degree $k$
supported on $R = nC \subset \overline \K_C$. As noted above, a semi-stable
vector bundle on $C$ of degree   $e = k + (n^2 -n)(1-g)$ determines a point
in $\text{Nilp}(\overline \H) $, i.e. there is natural inclusion $\U_C =
 \U_C(n, e) \subset \text{Nilp}(\overline \H)$; in the classical Hitchin
setting, this is just embedding of the zero-section in $\T^* \U_C$. We
write Nilp$(\H)$ for the corresponding locus of stable sheaves.  Similarly,
suppose the curve $C$ lies as a hyperplane section of a K3 surface $S$. We
define the nilpotent subvariety of Mukai space to be the set
$$\text{Nilp}(\overline \M) = \overline M^{-1}([R])$$
 consisting of all semi-stable sheaves of the appropriate numerical type
supported on the subscheme $R = nC \subset S$, where as above $\overline \M
= \overline \M^k_{|nC|}(S)$. The analogous locus of stable sheaves is
denoted  by Nilp$(\M) \subset \M$. Again one has an inclusion
$\U_C \subset \text{Nilp}(\overbar \M)$. Observe also that the nilpotent
subvariety Nilp$(\overbar \M)$ degenerates to a subvariety of Nilp$(\overbar
\H)$ under the deformation (1.6).

Our aim here is to prove the following:
\proclaim{\bf Theorem 2.1}  Assume that $k$ and $n$ are coprime. Then {\rm
Nilp}$(\M)$ is a Lagrangian subvariety of $\M$.
\endproclaim
\ni Recall from \cite{L}, Appendix A, that this means that Nilp$(\M)$
contains a Zariski-open dense set which is everywhere of half the
dimension of $\M$, and on which the symplectic form vanishes.

We start with
\proclaim{\bf Lemma 2.2} Under the numerical hypotheses of (2.1),
 every component of {\rm Nilp}$(\overline \M)$ has dimension $\tilde g$.
\endproclaim
\demo{Proof}
Let $\overline \H_0$ [resp. $\H_0$] denote the open subset of $\overline
\H$ [resp. $\H$] parametrizing semi-stable [resp. stable] sheaves supported
on $\K_C \subset \overline \K_C$. Note that sheaves on $\overline \K_C - \{
\text{vertex} \}$ are exactly those corresponding via the spectral
construction to Higgs pairs $(E, \phi : E \lra E\otimes K_C)$. We claim to
begin with that under the numerical hypotheses of the Theorem,
$$\overline \H_0 = \H_0. \tag *$$
In fact, suppose that
$\E$ is sheaf of pure dimension one on $\overline \K_C - \{ \text{vertex}
\}$. Then stability (or semi-stability) of $\E$ is in the first instance
defined
\`a la Gieseker by means of the Hilbert polynomial $p_{\E}(t) = \chi(
\E(tC))$ determined by the polarization
$C \subset \overline  \K_C$. However $C$ is linearly equivalent on $\K_C$ to
a divisor
$(2g-2)F$ where $F$ is the  pull-back of a divisor of
degree one on $C$ under the bundle projection $\K_C \lra C$, and it is
equivalent to compute stability with respect to the Hilbert polynomial
$$q_\E(t) = \chi(\E(t  F)).$$
Now
$$q_\E(t) = nt + \left ( k + n^2(1-g)
\right ),$$ and hence if $k$ and $n$ are relatively prime, then so are the
two coefficients $r = n$ and $a = k + n^2(1-g)$ of $q_\E(t)$. But  from the
 definition of $q$-stability it follows immediately that if
$\E$ is semi-stable, then it  is automatically stable. [In effect, we are
proving that (semi-)stability of $\E$ is equivalent to (semi-)stability of
the corresponding Higgs pair $(E, \phi)$.]

We next wish to invoke Laumon's theorem that the Hitchin nilpotent
cone is Lagrangian. The main theorem of [L] states that the nilpotent cone is
Lagrangian in the moduli stack $T^*\text{Fib}$, which  parametrizes
isomorphism classes of arbitrary Higgs pairs $(E , \phi)$. In particular,
each component of this nilpotent cone has dimension (as a stack) $\tilde g -
1$, which is the dimension of Fib. (The difference of
$1$ is due to the presence of scalar automorphisms on every vector bundle.)
But as we have just seen, points of $\overline \H_0$ are given by
{\it stable} pairs $(E,\phi)$, and since stable pairs are simple by
\cite{S},\S1,  we conclude that each component of $\nilp(\overline \H)$ is
Lagrangian of dimension $\tilde g$ in ${\overline \H_0}$.

 Now under the deformation  (1.6), Nilp$(\overline \M)$ specializes to
Nilp$(\overline \H)$. Hence it follows from the semi-continuity of fibre
dimensions that every component of Nilp$(\overline \M)$ has dimension $\le
\tilde g$. On the other hand, Lemma 1.3 shows that Nilp$(\overline \M)$ is
itself a specialization of the Lagrangian subvarieties $\overline
\M^{-1}([D])$ for general $D \in |nC|$, and this implies the reverse
inequality. \qed \enddemo

\ni{\bf Remark 2.3.} It seems plausible that Lemma 2.2 remains valid
without the numerical hypotheses of the Theorem.  However we
do not know how to rule out the possibility that there are whole components
of $\nilp(\overline \H)$ contained in the singular locus of $\overline \H$,
in which case  Laumon's theorem does not seem to apply. Note in any event
that it is only in Lemma 2.2 that we use that $k$ and $n$ are coprime.

 Keeping the hypotheses of the Theorem, it remains to show that the
symplectic form on $\M$ vanishes on any component $F$ of Nilp$(\M)$. Denote
by $\overline F$ the closure of
$F$ in Nilp$(\overline \M)$, so that $\overline F$ is an irreducible
component of $\nilp(\overline \M)$. Now $\M$ has dimension
$2 \tilde g$ at every point of $F$, so it follows from Lemma 2.2 that
$\overline F$ is contained in an irreducible component $W$ of $\overline
\M$ that maps onto $\overline B  = |nC|$. Choose a smooth
curve $T \subset \overline B$ passing through
$[R] = [nC] \in \overline B$, and meeting the open subset of $\overline B$
parametrizing smooth members of $|nC|$. Let $W_T =  W \cap \overline
M^{-1}(T) \subset \overline \M$, and denote by $m : W_T \lra T$ the
projection. Then every component of $W_T$ has dimension $\ge \tilde g + 1$,
whereas $\overline F$, which has dimension $\tilde g$, is an irreducible
component of $  \nilp(\overline \M) \supseteq m^{-1}([R])$. It follows that
there is an irreducible component $Z$ of $W_T$ containing $\overline F$ and
mapping onto $T$. Write $g : Z \lra T$ for the projection. Thus $\overline F$
is a component of $g^{-1}(0)$, whereas for $t$ in a punctured
neighborhood of $[R] \in T$, it follows from (1.3) that the fibre $Z_t =
g^{-1}(t)$ is a Lagrangian submanifold of $\M$.  Therefore the assertion
follows from the following lemma, which states in effect that a limit of
Lagrangian submanifolds is a Lagrangian subvariety.

\proclaim{\bf Lemma 2.4}  Let $(M, \omega)$ be a symplectic variety, let $T$
be (the germ of) a smooth curve with a marked point $0 \in T$. Suppose that
$Z \subset M \times T$ is an irreducible variety such that for $0 \ne t \in
T$ the fibre $Z_t \subset M$ is a smooth subvariety on which the symplectic
form $\omega$  vanishes. Then $\omega$ also vanishes at the general point of
any component of the special fibre $Z_0 \subset M$. \endproclaim
\demo{Proof}  Let $p : Z \lra T$ denote projection onto the second factor.
By Mumford's semi-stable reduction theorem \cite{K}, after a base change
$(T^\prime , 0^\prime) \lra (T, 0)$ and blowings up over the central fibre,
we can construct a new family $q : Y \lra T^\prime$ such that the central
fibre $q^{-1}(0^\prime)$ is a reduced normal crossing divisor, and all other
fibres are smooth. Then the sheaf
$\Omega^2_{Y/T^\prime}$ of relative two-forms along $q$ is locally free
outside the subset $G \subset Y$ of codimension $\ge 2$ where $q$ fails to be
smooth.  Let $f : Y \lra M$ denote the composition $Y \lra Z \lra M$. Then $
\eta =_{\text{def}} f^*(\omega)$ vanishes on the general fibre of $q$, and
hence vanishes off $G$. In particular, $\eta$
vanishes at the general point of each component $D_i$ of $q^{-1}(0^\prime)$.
But each component $W$ of $Z_0 \subset Z$ is dominated by some $D_i$, so
$\omega$ must vanish at the general point of $W$. \qed \enddemo

This completes the proof of Theorem 2.1.

\bl
\ni {\bf \S3. Comparison of Nilpotent Cones in Rank Two.} In this section,
we state without proof a detailed description of the nilpotent cones in the
case $n = 2$. For simplicity we limit attention to stable sheaves, and the
classical Hitchin setting, although we do not need the numerical hypotheses
of Theorem 2.1.

 We start with the Hitchin system. Consider the first infinitesimal
neighborhood $R = 2C \subset \K_C$. Via the spectral construction one has an
identification
$$
\text{Nilp}( \H) =  H^{-1}(0) = \left \{ \  (E, \phi) \
\bigg |
\gathered
\text{E a stable rank 2 bundle  on $C$,}\\
\phi : E \lra E \otimes K_C \text{ with } \phi^2 = 0 \
\endgathered  \ \  \right
\}.$$ Laumon \cite{L} has enumerated the components of the nilpotent cone in
general, but in the present case the description is elementary to obtain by
hand. In fact, suppose
$0 \ne \phi$ is nilpotent. Then $\phi$ has generic rank one.  Let $D \subset
C$ be the effective divisor on which $\phi$ vanishes, and denote by $A$ the
line bundle $\text{im} \phi \subset E \otimes K_C$. Then
$E$ sits in the exact sequence $$
0 \lra A(D) \otimes K_C^{-1} \overset{\alpha}
\to \lra E \lra A \lra 0,$$ and
$\phi$ factors as the composition
$$
\CD E @>\alpha >> A @. \\ @. @VV{\cdot D }V @. \\ @. A(D) @>{\alpha \otimes
1}>>  E \otimes K_C.
\endCD
$$
 Such nilpotent endomorphisms may be parametrized as follows. Set $e =
\text{deg }E$,  and fix an integer
$$0 \le d < 2g - 2, \ \ \text{with} \ \ d \equiv e \pmod{2}.$$ Put $a = (e +
2g - 2 - d)/2$, and define:
$$
\N_d = \big \{ (A, D, \epsilon )
 \bigm| A \in Pic^a(C), \
		D \in  Sym^d(C), \
\epsilon \in \text{\rm Ext}^1(A, A(D) \otimes K_C^{-1})
\big \}.
$$ Note that there is a natural map
$$\pi : \N_d \lra Pic^a(C) \times Sym^d(C)$$ which realizes $\N_d$ as the
total space of a vector bundle whose fibre over the point $(A, D)$ is the
vector space
$H^1(\O(D) \otimes K_C^{-1})$. In the special case at hand, we may state one
of Laumon's results as follows:
\proclaim{Proposition 3.1}  The irreducible components of the nilpotent
cone {\rm Nilp}$  \H$  consist of the moduli space $\U_C(2,e)$
of stable bundles, together with the  Zariski closures in
$ \H$ of  suitable Zariski-open subsets
$N^+_d \subset \N_d$. \qed \endproclaim
\ni (The open subsets $N_d^+$ arise because not every rank two bundle
determined by an element in $\N_d$ is stable.)

Turning to the Mukai setting, recall that the nilpotent subvariety
Nilp$( \M)$ para- metrizes stable sheaves of numerical rank one
and degree $k$ on the infinitesimal neighborhood $R = 2C \subset S$ of $C$
in the K3 surface $S$. At this point the question no longer has anything to
do with K3 surfaces, and it clarifies matters to generalize slightly.

	Let $C$ be a smooth curve of genus $g$, and let  $R$ be any ``ribbon" on
$C$ with normal bundle $K_C$, i.e. the scheme arising from a double structure
on $C$, with normal bundle $K_C$. Thus $\O_R$ sits in an exact sequence:
$$ 0 \lra K_C^{-1} \lra  \O_R \lra \O_C \lra 0.$$ Such double structures are
classified by $H^1(C, K_C^{-2})$, and given $r \in H^1(C, K_C^{-2})$ we write
$R_r$ for the corresponding scheme. For example, $R_0$ is just
the first infinitesimal neighborhood of the zero-section $C \subset
\K_C$. Fixing some polarization on $R$, we define $  \M(r) =
  \M^k(R_r)$ to be the moduli space of stable (in particular,
Cohen-Macaulay) sheaves of numerical rank one and degree $k$ on $R_r$.

	The sheaves in question have one of two types. First, a rank two
stable vector bundle on $C$ of degree $e = k + 2 - 2g$, considered as an
$\O_R$ module, has numerical rank one and degree $k$. The second type of
sheaves consists of stable $\O_R$-modules $\E$ whose restrictions $\E
\otimes \O_C$ to $C \subset R$ have rank 1. Fix such a sheaf $\E$. One shows
that canonically associated to $\E$ there is an effective divisor $D \subset
C$ on $C$ supported on the set where $\E$ fails to be locally isomorphic to
$\O_R$:  we write $D = \text{Sing}(\E)$. Put  $$A = \E \otimes \O_C \ /
\ \text{torsion}.$$ One has an exact sequence
$0 \lra A(D) \otimes K_C^{-1} \lra \E \lra A \lra 0$
of $\O_R$-modules. Set
$$\N_d(r) = \big \{ \E \in \M(r) \bigm | \text{ deg Sing}(E) = d \big \}.$$
Then there is a morphism
$$\pi : \N_d(r) \lra Pic^a(C) \times Sym^d(C) \quad \text{via} \quad \E
\mapsto (A, D),$$
where as above $D = \text{Sing}(\E)$, and $a = (k -d)/2$.

We may now state the result:
\proclaim{\bf Theorem 3.2} (i). The irreducible components of
$ \M^k(r)$ are the Zariski closures of
$$\N_{\infty}(r) =_{\text{def}} \U_C(2, e) \qquad (\text{for
$e = k + 2 - 2g$}) $$ and the sets $$\N_d(r)\qquad (\text{for $0 \le d < 2g -
2$ and
$d
\equiv k \pmod{2}$}).$$

\ni (ii). The map $\pi : \N_d(r) \lra Pic^a(C) \times Sym^d(C)$ is an affine
bundle, with underlying vector bundle $\N_d = \N_d(0)$. In general it
is not itself a vector bundle.
\endproclaim
\ni For example consider the case $d = 0$. Then $\N_0(r) = Pic^k(R_r)$, and
the morphism $\pi : Pic^k(R_r) \lra Pic^a(C)$ comes  from the exact sequence
$$
0 \lra H^1(C, K_C^*) \lra Pic^k(R_r) \lra Pic^a(C) \lra 0 \tag 3.3
$$
of algebraic groups. The last statement of the Theorem is illustrated by
the fact that $Pic^k(R_r)$  has the structure of an affine bundle, but not
in general a vector bundle, over $Pic^a(C)$. Note that the
Theorem can be seen simply as a statement about the geometry of double
curves, but it seems to us that it is only in the context of the
Hitchin-Laumon nilpotent cone that it becomes natural. The proof of the
Theorem is not terribly difficult but it is fairly long, and we will not
give details here.  We remark that a similar analysis holds for the moduli
space of stable sheaves on a ``ribbon" $R$ with any normal bundle $N$, at
least when $N$ has sufficiently large degree.

There is some interesting additional geometry connected with this situation.
Let $\E$ be a sheaf on $R$, corresponding to a point in
$\N_d(r)$, and set $D = \text{Sing} \ \E$. One can show that there exists a
line bundle $\L$ on $R$ such that
$$ \E = \text{elm}_D(\L) =_{\text{def}} \ker \{ \L \lra \L \otimes \O_D
\}.$$
This globalizes to a morphism
$$f : Pic^k(R) \times Sym^d(C) \lra \N_d(r)$$
given  by $f( \L, D) = \text{elm}_D(\L)$. In fact, $f$ is a map of affine
bundles over $Pic^a(C) \times Sym^d(C)$, and for some questions it reduces
the analysis of $\N_d(r)$ to the case $d = 0$.

Recall next that if $X$ is any variety, and if $E$
is a vector bundle on $X$, then the set of all isomorphism classes of affine
bundles on $X$ with underlying linear structure $E$ is classified by the
cohomology group $H^1(X, E)$. It is an amusing exercise to compute the
characteristic classes of the various affine bundles appearing in the
Theorem. When $d = 0$, for example,  the bundle (3.3) is classified by an
element in
$$
H^1(Pic(C), H^1(K_C^{-1}) \otimes_\C \O_{Pic(C)}) = H^1(C, K_C^{-1})
\otimes H^1(C, \O_C).$$
As one might expect, the element in question is just the image of the
``ribbon class" $r$ under the natural map $H^1(K_C^{-2}) \lra H^1(\O_C)
\otimes H^1(K_C^{-1})$ dual to the multiplication $H^0(K_C) \otimes
H^0(K_C^2) \lra H^0(K_C^3)$.

Finally, note that it follows immediately from (3.2.ii) that every component
of $\M(r)$ has dimension $g + d + \text{rank}\N_d = 4g - 3 = \tilde g$.
Returning to the Mukai nilpotent cone $\nilp(\M)$ determined by a
$K3$-``ribbon" $R = 2C \subset S$, this implies that Lemma 2.2 remains
valid when $n = 2$ without the hypothesis that $k$ be odd. Hence in the
rank two case, we do not need any restrictions on $k$ in Theorem 2.1.

\bl
\ni{\bf \S 4. Concluding Remarks and Open Questions. } It would be quite
interesting to generalize the analysis of \S 3 to higher rank. Let $R$ be a
``ribbon" (or ``tape") of order $n$ and normal bundle $N$, i.e. a
multiplicity $n$ structure on a smooth curve $C$ which looks locally like the
$(n-1)^{\text{st}}$ infinitesimal neighborhood of $C$ on a surface. Thus we
suppose that $R$ is filtered by subschemes
$$C = C_{1} \subset C_2 \subset \dots \subset C_{n-1} \subset C_n = R,$$
where the $C_i$ sit in exact sequences
$$ 0 \lra N^{-(i-1)} \lra \O_{C_i} \lra \O_{C_{i-1}} \lra 0.$$
For example, one has the ``split" ribbon $R_0$, i.e. the $(n-1)^{\text{st}}$
neighborhood of the zero-section in the normal bundle $\N$.   Consider a
Cohen-Macaulay $\O_R$-module $\E$ of numerical rank one. Define
$$\nu_i = \text{rank}_{\O_C} \ker \{ \E \otimes \O_{C_i} \lra \E \otimes
\O_{C_{i-1}} \}$$
(and $\nu_1 = \text{rank}(\E \otimes \O_C)$).
One has $\nu_1 \ge \nu_2 \ge \dots \ge \nu_n$, and we call the
vector $\nu = (\nu_1, \dots ,\nu_n)$ the {\it type} of $\E$. For example,
when $n = 2$ there are two possible types: $(2,0)$, corresponding to a rank
two vector bundle on $C$, and $(1,1)$, corresponding to a sheaf $\E$ whose
restriction to $C$ has generic rank one. As another example, on the split
ribbon $R_0$,
$\E$ is given by a vector bundle $E$ of rank $n$ on $C$  together with a
nilpotent endomorphism $\phi : E \lra E \otimes N$, and then
$\nu_i = \text{rank}( \ker \phi^i / \ker \phi^{i-1})$. One can associate to
a sheaf $\E$ of type $\nu$ vector bundles $E_i$ on $C$ of rank $\nu_i$
(generalizing the line bundle $A$ in \S 3) and sky-scraper sheaves
$\Delta_i$ on $C$ (generalizing the divisor $D = \text{ Sing}(E)$ in \S3).
The type -- and, within a given type, the degrees of the $E_i$ and
$\Delta_i$ -- give discrete invariants of the components of the space of
stable sheaves on $R$. For type $(1,1, \dots, 1)$ one can work out the
picture  in some detail, much as in \S 3. However for general types, the
analysis is less clear.  One can hope that the sort of construction with
``elementary transformations" indicated at the end of the previous
section can at least reduce the study of the general case to the
``defectless" situation where all $\Delta_i = 0$.

  Much more generally, one can attempt to replace the $K3$ surface by
a higher dimensional symplectic variety  $S$. It is not clear at the moment
exactly when one should expect to have an analogue of Mukai's theorem,
 i.e. a natural symplectic structure on appropriate moduli spaces of sheaves
on $S$.  The case of vector bundles on $S$ was proved by Kobayashi \cite{Ko}.
There are however examples of such moduli spaces (parametrizing line bundles
on a divisor in $S$) which have odd dimension, and so cannot admit any
algebraic symplectic structure, cf. \cite{DM2}.  A beautiful new idea of
Markman \cite{Mkm} is that Mukai's results do extend to the moduli of
Lagrangian sheaves, i.e. those sheaves on a symplectic variety $S$
whose support in $S$ is itself  Lagrangian. Further, a version of this space
comes with a support map to a Lagrangian-Hilbert base $B$, and this map makes
it into an  algebraically completely integrable system. More precisely,  there
is a
closed two-form on the non-singular locus of the moduli space of Lagrangian
sheaves, which can be constructed in terms of a field of cubics on the base $B$
as in Remark 1.4, cf. \cite{DM1}.  Under some mild conditions, this form is
non-degenerate at all points which correspond to line bundles over a
non-singular
Lagrangian support. There are several variations which replace the symplectic
structures (on $S$ and on the resulting moduli spaces) with Poisson or
quasisymplectic structures, respectively.  This whole picture can be considered
as a non-linear analogue, and sometimes as a deformation, of Simpson's moduli
space of Higgs bundles on a variety \cite{S}, or of the more general moduli
spaces
of vector bundles with arbitrarily twisted endomorphisms studied e.g. in
\cite{DM2},
in the same sense as Mukai's space appears here as a deformation of Hitchin's.

\bl
\ni{\bf References.}

\vskip 20pt
\bl
\settabs\+University of Illinois at Chicago and now is the time  \cr
\+ Ron DONAGI \cr
\+ Department of Mathematics \cr
\+ University of Pennsylvania \cr
\+ Philadelphia, PA  19104 \cr
\+ e-mail:  donagi{\@}math.upenn.edu \cr

\vskip 7pt

\+ Lawrence EIN \cr
\+ Department of Mathematics \cr
\+ University of Illinois at Chicago \cr
\+ Chicago, IL  60680 \cr
\+ e-mail: U22425$\%$UICVM.BITNET \cr

\vskip 7 pt

\+ Robert LAZARSFELD \cr
\+ Department of Mathematics \cr
\+ University of California, Los Angeles \cr
\+ Los Angeles, CA  90024 \cr
\+ e-mail: rkl$\@$math.ucla.edu \cr

\end